\documentclass[11pt, a4paper]{article}
\usepackage{amsmath,amssymb,latexsym,color}
\usepackage{authblk}
\usepackage[mathscr]{eucal}
\usepackage[colorlinks,
            linkcolor=blue,
            anchorcolor=green,
            citecolor=magenta
           ]{hyperref}
\newtheorem{thm}{Theorem}[section]

\newtheorem{lemma}[thm]{Lemma}

\newtheorem{ass}{Assumption}
\newtheorem{remark}{Remark}

\newtheorem{example}{Example}[section]
\newtheorem{defin}{Definition}[section]

\newcommand{\proof}{{\it Proof.\quad}}
\newcommand{\qed}{\hfill\Box\medskip}

\usepackage{CJK}

\begin{document}
%\begin{CJK*}{GBK}{song}

\renewcommand{\baselinestretch}{1.3}
%%%%%%%%%%%%%%%%%%%%%%%%%%%%%%%%%%%%%%%%%%%%%%%%%%%%%%%%%%%%%%%%%%%%%%%%%%%%%%%%%%%%%%%%
%%%%%%%%%%%%%%%%%%%%%%%%%%%%%%%%%%%%%%%%%%%%%%%%%%%%%%%%%%%%%%%%%%%%%%%%%%%%%%%%%%%%%%%%
\title{\bf The structure of large non-trivial $t$-intersecting families for finite sets}
\author[1,2]{Mengyu Cao\thanks{E-mail: \texttt{caomengyu@mail.bnu.edu.cn}}}
\author[1]{Benjian Lv\thanks{Corresponding author. E-mail: \texttt{bjlv@bnu.edu.cn}}}
\author[1]{Kaishun Wang\thanks{E-mail: \texttt{wangks@bnu.edu.cn}}}
\affil[1]{\small Laboratory of Mathematics and Complex Systems (Ministry of Education), School of Mathematical Sciences, Beijing Normal University, Beijing 100875, China}
\affil[2]{\small Department of Mathematical Sciences, Tsinghua University, Beijing 100084, China}
 \date{}
 \maketitle

\begin{abstract}
In this paper, we describe the structure of maximal non-trivial uniform  $t$-intersecting families with large size for finite sets. In the special case when $t=1$, our result gives rise to Kostochka and Mubayi's result in 2017.

\medskip
\noindent {\em AMS classification:} 05D05, 05A10

\noindent {\em Key words:} Erd\H{o}s-Ko-Rado Theorem; Hilton-Milner Theorem;  $t$-intersecting family; $t$-covering number

\end{abstract}

\section{Introduction}
Let $n$ and $k$ be integers with $1\leq k\leq n.$ Write $[n]=\{1,2,\ldots,n\}$ and denote by ${[n]\choose k}$ the family of all $k$-subsets of $[n].$ For any positive integer $t$, a family $\mathcal{F}\subseteq {[n]\choose k}$ is said to be \emph{$t$-intersecting} if $|A \cap B|\geq t$ for all $A, B\in\mathcal{F}.$ A $t$-intersecting family is called \emph{trivial} if all its members contain a common specified $t$-subset of $[n]$, and \emph{non-trivial} otherwise.

The famous Erd\H{o}s-Ko-Rado Theorem gives the maximum size of a $t$-intersecting family and shows further that any $t$-intersecting family with maximum size is a trivial family consisting of all $k$-subsets that contain a fixed $t$-subset of $[n]$ for $n>n_0(k,t)$ \cite{Erdos-Ko-Rado-1961-313}. It is known that the smallest possible such function $n_0(k, t)$ is $(t +1)(k-t +1).$ This was proved by Frankl \cite{Frankl-1978} for $t\geq 15$ and subsequently determined by Wilson \cite{Wilson-1984} for all $t$. In \cite{Frankl-1978}, Frankl also made a conjecture on the maximum size of a $t$-intersecting family of $k$-subsets of $[n]$ for any positive integers $t,k$ and $n$. This conjecture was partially proved by Frankl and F\"{u}redi in  \cite{Frankl--Furedi-1991} and completely settled by Ahlswede and Khachatrian in \cite{Ahlswede-Khachatrian-1997}.

Determining the structure of non-trivial $t$-intersecting families of $k$-subsets of $[n]$ with maximum size was a long-standing problem. The first result is the Hilton-Milner Theorem \cite{Hilton-Milner-1967} which describes the structure of such families for $t=1$. A significant step was taken in \cite{Frankl-1978-1} by Frankl, who determined such families for $t\geq 2$ and $n>n_1(k,t)$. In \cite{Frankl-Furedi-1986}, Frankl and F\"{u}redi gave a short and elegant proof for the Hilton-Milner Theorem by using the shifting technique, and also asked whether $n_1(k,t)<ckt$ holds.  Ahlswede and Khachatrian \cite{Ahlswede-Khachatrian-1996} answered this question and gave a complete result on non-trivial intersection problems for finite sets.

Recently, other maximal non-trivial $1$-intersecting families with large size have been studied. In \cite{Han-Kohayakawa}, Han and Kohayakawa determined the structure of the third largest maximal $1$-intersecting families of $k$-subsets of $[n]$ with $6\leq 2k<n$. They also mentioned that it would be natural to investigate such problem for $t$-intersecting families. In \cite{Kostochka-Mubayi}, Kostochka and Mubayi  described the structure of $1$-intersecting families of $k$-subsets of $[n]$ for large $n$ whose size is quite a bit smaller than the bound ${n-1\choose k-1}$ given by the Erd\H{o}s-Ko-Rado Theorem.

Define the $t$-\emph{covering number} $\tau_t(\mathcal{F})$ of a family $\mathcal{F}\subseteq {[n]\choose k}$ to be the minimum size of a subset $T$ of $[n]$ such that $|T\cap F|\geq t$ for any $F\in \mathcal{F}$. Let $\mathcal{F}\subseteq{[n]\choose k}$ be any $t$-intersecting family. Note that $t\leq \tau_t(\mathcal{F})\leq k$, and $\mathcal{F}$ is trivial if $\tau_t(\mathcal{F})=t.$ In \cite{Frankl-1978-1}, in order to describe the structure of non-trivial $t$-intersecting families with maximum size for $t\geq 2$, Frankl defined the family $\mathcal{F}^{(t)}$ and the base $\mathcal{B}$ of $\mathcal{F}.$ It is straightforward to verify that the $t$-covering number of a maximal non-trivial $t$-intersecting family equals to $\ell_1$ in the decomposition of its base $\mathcal{B}.$  Observe that when $t=1$ that is the well-known covering number of an intersecting family. We refer the readers to \cite{Frankl-1980,Frankl--Ota--Tokushige-1991,Furedi-1988,Furuya--Takatou-1988} for more results about the covering number.

In this paper, we consider maximal non-trivial $t$-intersecting families with large size for any positive integer $t$. If $t=k-1$, it is well known that any maximal non-trivial $(k-1)$-intersecting family is a collection of all $k$-subsets containing a fixed $(k-1)$-subset or a collection of all $k$-subsets contained in a fixed $(k+1)$-subset. Thus, we only consider the case with $1\leq t\leq k-2$. To present our result let us first introduce the following two constructions of $t$-intersecting families of $k$-subsets of $[n]$.

\

\noindent\textbf{Family I.}\quad Let $X$, $M$ and $C$ be three subsets of $[n]$ such that $X\subseteq M\subseteq C$, $|X|=t$, $|M|=k$ and $|C|=c$, where $c\in\{k+1,k+2,\ldots,2k-t,n\}$. Denote
\begin{align*}
\mathcal{H}_1(X,M,C)=\mathcal{A}(X,M)\cup\mathcal{B}(X,M,C)\cup\mathcal{C}(X,M,C),
\end{align*}
where
\begin{align*}
\mathcal{A}(X,M)=&\left\{F\in{[n]\choose k}\mid X\subseteq F,\ |F\cap M|\geq t+1\right\},\\
\mathcal{B}(X,M,C)=&\left\{F\in{[n]\choose k}\mid F\cap M=X,\ |F\cap C|=c-k+t\right\},\\
\mathcal{C}(X,M,C)=&\left\{F\in{C\choose k}\mid |F\cap X|=t-1,\ |F\cap M|=k-1\right\}.
\end{align*}

\

\noindent\textbf{Family II.}\quad Let $Z$ be a $(t+2)$-subset of $[n]$. Define
\begin{align*}
\mathcal{H}_2(Z)=\left\{F\in{[n]\choose k}\mid|F\cap Z|\geq t+1\right\}.
\end{align*}

It is straightforward to verify that Families I and II are $t$-intersecting families with $t$-covering number $t+1$. Observe that the size of each this family only depends on $|X|$, $|M|$, $|C|$ and $|Z|$. Let  $h_1(t,k,c)=|\mathcal{H}_1(X,M,C)|$, where $c=|C|\in\{k+1,k+2,\ldots,2k-t,n\}$;  $h_2(t+2)=|\mathcal{H}_2(Z)|$.

\begin{remark}\label{s-obs1} {\rm Suppose $X$, $M$ and $C$ are three subsets of $[n]$ satisfying the condition in Family I. If $|C|=k+1$, then
$$
\mathcal{H}_1(X,M,C)=\left\{F\in{[n]\choose k}\mid X\subseteq F,\ |F\cap C|\geq t+1\right\}\cup {C\choose k}.
$$
 If $t=k-2$, then $\mathcal{H}_1(X,M,[n])=\mathcal{H}_2(M)$ and $h_1(k-2,k,n)=h_2(k)$. }
\end{remark}

Our main result describes the structure of all maximal non-trivial uniform $t$-intersecting  families with large size for finite sets.
\begin{thm}\label{s-main-1}
Let $1\leq t\leq k-2$, and $\max\left\{{t+2\choose 2},\frac{k-t+2}{2}\right\}\cdot(k-t+1)^2+t\leq n$. If  $\mathcal{F}\subseteq{[n]\choose k}$ is a maximal non-trivial $t$-intersecting family with
$$
|\mathcal{F}|\geq (k-t){n-t-1\choose k-t-1}-{k-t\choose 2}{n-t-2\choose k-t-2},
$$
then one of the following holds.
\begin{itemize}
\item[{\rm(i)}] $\mathcal{F}=\mathcal{H}_1(X,M,C)$ for some $t$-subset $X$, $k$-subset $M$ and $c$-subset $C$ of $[n]$, where $c\in\{k+1,k+2,\ldots,2k-t,n\}$.
\item[{\rm(ii)}] $\mathcal{F}=\mathcal{H}_2(Z)$ for some $(t+2)$-subset $Z$ of $[n],$ and $\frac{k}{2}-1\leq t\leq k-2$.
 \end{itemize}
\end{thm}

In the special case when $t=1$, Theorem~\ref{s-main-1} gives rise to Theorem 7 in \cite{Kostochka-Mubayi}. From Lemmas~\ref{s-lem3}, \ref{s-lem6}, \ref{s-lem2} and \ref{s-lem2-2}, one can determine the sequence of these families according to their sizes.

The rest of this paper is organized as follows. In the next section we will give some properties of the maximal $t$-intersecting families with $t$-covering number $t+1$, and prove a number of inequalities for the sizes of Families I and II. In Section 3 we will prove some upper bounds for the sizes of non-trivial $t$-intersecting families using their $t$-covering number. After these preparations we will prove Theorem~\ref{s-main-1} in Section 4.

\section{$t$-intersecting families with $t$-covering number $t+1$}
In this section, we will give some properties of the maximal $t$-intersecting families with $t$-covering number $t+1$ and prove a number of inequalities for the sizes of Families I and II.
\begin{ass}\label{s-hyp1}
	Let $1\leq t\leq k-2$ and $2k\leq n$, let $\mathcal{F}\subseteq {[n]\choose k}$ be a maximal $t$-intersecting family with $\tau_t(\mathcal{F})=t+1.$ Define
	$$
	\mathcal{T}=\left\{T\in{[n]\choose t+1}\mid |T\cap F|\geq t\ for\ any\ F\in\mathcal{F}\right\}.
	$$
\end{ass}

Observe that the $\mathcal{T}$ in Assumption~\ref{s-hyp1} is the family $\mathcal{B}_{\ell_1}$ with $\ell_1=t+1$ in the decomposition of the base $\mathcal{B}$ of $\mathcal{F}$ in \cite{Frankl-1978-1}. Applying the method in the proof of Theorem~1 in  \cite{Frankl-1978-1}, we firstly give several properties of the maximal $t$-intersecting families with $t$-covering number $t+1$.
\begin{lemma}\label{s-lem3-1}
	Let $n,\ k,\ t,\ \mathcal{F}$ and $\mathcal{T}$ be as in Assumption~\ref{s-hyp1}. Then $\mathcal{T}$ is a $t$-intersecting family with $t\leq\tau_t(\mathcal{T})\leq t+1$. Moreover, the following hold.
	\begin{itemize}
		\item[{\rm (i)}] If $\tau_t(\mathcal{T})=t$, then there exist a $t$-subset $X$ and an $l$-subset $M$ of $[n]$ with $X\subseteq M$ and $t+1\leq l\leq k+1$ such that \begin{align}\label{s-equ1-1}
		\mathcal{T}=\left\{T\in{M\choose t+1}\mid X\subseteq T\right\}.
		\end{align}
		\item[{\rm (ii)}] If $\tau_t(\mathcal{T})=t+1$, then there exists a $(t+2)$-subset $Z$ of $[n]$ such that $\mathcal{T}={Z\choose t+1}.$
	\end{itemize}
\end{lemma}
\proof For any $T\in\mathcal{T}$, by maximality of $\mathcal{F}$, $\mathcal{F}$ contains all $k$-subsets of $[n]$ containing $T$. For any $T_1,T_2\in \mathcal{T},$ if $|T_1\cap T_2|<t,$ then there must exist $F_1,F_2\in\mathcal{F}$ such that $T_1\subseteq F_1$, $T_2\subseteq F_2$ and $|F_1\cap F_2|<t$ from $2k\leq n.$ That is impossible as $\mathcal{F}$ is maximal $t$-intersecting. Hence $|T_1\cap T_2|\geq t,$ and $\mathcal{T}\subseteq {[n]\choose t+1}$ is a $t$-intersecting family with $t\leq\tau_t(\mathcal{T})\leq t+1.$
\
\medskip
\

(i)\quad Suppose that $\tau_t(\mathcal{T})=t.$ Then there exists a $t$-subset $X$ of $[n]$ such that $X$ is contained in every $(t+1)$-subset in $\mathcal{T}.$ Assume that $M=\cup_{T\in\mathcal{T}}T$ and $|M|=l.$ It suffices to prove (\ref{s-equ1-1}) and $t+1\leq|M|\leq k+1$. Since $\tau_t(\mathcal{F})=t+1$, we have $\mathcal{F}\setminus\mathcal{F}_X\neq\emptyset.$ Let $F^\prime$ be any $k$-subset in $\mathcal{F}\setminus\mathcal{F}_X.$ Observe that
$|X\cap F^\prime|\leq t-1.$ For any $T\in\mathcal{T},$ since $X\subseteq T$ and $|T\cap F^\prime|\geq t,$ we have $|X\cap F^\prime|=t-1$ and $|T\cap (X\cup F^\prime)|\geq t+1,$ which imply that $|X\cup F^\prime|=k+1$ and $T\subseteq X\cup F^\prime.$ Hence $M=\cup_{T\in\mathcal{T}}T\subseteq X\cup F^\prime$ and $t+1\leq l\leq k+1.$ It is clear that $\mathcal{T}\subseteq\{T\in{M\choose t+1}\mid X\subseteq T\}.$ Let $T^\prime$ be any $(t+1)$-subset of $M$ with $X\subseteq T^\prime.$ For any $F\in\mathcal{F}$, if $X\subseteq F,$ then $|T^\prime\cap F|\geq t$; if $X\nsubseteq F,$ by above discussion, then $T^\prime\subseteq X\cup F,$ which implies that $|T^\prime\cap F|\geq t$ from $|X\cup F|=k+1.$ Hence $T^\prime\in\mathcal{T}$ and (\ref{s-equ1-1}) is proved.
\
\medskip
\

(ii)\quad Suppose that $\tau_t(\mathcal{T})=t+1$. Let $A,B,C\in\mathcal{T}$ be distinct subsets such that $A\cap B$, $A\cap C$ and $B\cap C$ are pair-wise distinct. Since $\mathcal{T}$ is $t$-intersecting, we have $|A\cap B|=|A\cap C|=|B\cap C|=t,$ which implies that $C=(A\cap C)\cup(B\cap C)\subseteq A\cup B$ from $|C|=t+1.$ Hence, we get $A\cup C\subseteq A\cup B$ and $B\cup C\subseteq A\cup B$, which imply that $A\cup B=A\cup C=B\cup C.$

Since $\tau_t(\mathcal{T})=t+1$, there exist three distinct subsets $T_1,T_2,T_3\in\mathcal{T}$ such that $T_1\cap T_2$, $T_1\cap T_3$ and $T_2\cap T_3$ are pair-wise distinct.  For any $T\in\mathcal{T}\setminus\{T_1,T_2,T_3\},$ if $T\cap T_1=T\cap T_2=T\cap T_3,$ then $|T\cap T_1|=t,$ $T\cap T_1\subseteq T_2$ and $T\cap T_1\subseteq T_3,$ which imply that $T\cap T_1=T_1\cap T_2=T_1\cap T_3,$ a contradiction. Hence, there exist $T_i, T_j\in\{T_1,T_2,T_3\}$ such that $T\cap T_i\neq T\cap T_j,$ and $
T= (T\cap T_i)\cup (T\cap T_j)\subseteq T_1\cup T_2=T_1\cup T_3=T_2\cup T_3.$

Let $Z=T_1\cup T_2.$ Then $\mathcal{T}\subseteq{Z\choose t+1}.$ In the following, we show that ${Z\choose t+1}\subseteq\mathcal{T}.$ For any $F\in \mathcal{F},$ if $F\cap T_1=F\cap T_2=F\cap T_3,$ then $F\cap T_1\subseteq T_i$ for any $i\in\{1,2,3\}.$ That is impossible as $T_1\cap T_2,$ $T_1\cap T_3$ and $T_2\cap T_3$ are pair-wise distinct and $|F\cap T_1|\geq t$. Hence there exist $T_i,T_j\in\{T_1,T_2,T_3\}$ such that $F\cap T_i\neq F\cap T_j,$ implying that $|F\cap Z|\geq t+1.$ So for any $F\in\mathcal{F}$ and $T^\prime\in{Z\choose t+1},$ we have $|F\cap T^\prime|\geq t.$ Therefore, we have $\mathcal{T}={Z\choose t+1}$ as desired. $\qed$
\begin{lemma}\label{s-prop3}
	Let $n,\ k,\ t,\ \mathcal{F}$ and $\mathcal{T}$ be as in Assumption~\ref{s-hyp1}, and set $M=\cup_{T\in \mathcal{T}}T$. Suppose that $\tau_t(\mathcal{T})=t$, $|M|=k+1$ and $X$ is a $t$-subset of $[n]$ which is contained in each $T\in\mathcal{T}.$ Then $\mathcal{F}=\{F\in {[n]\choose k}\mid X\subseteq F,\ |F\cap M|\geq t+1\}\cup{M\choose k}.$
\end{lemma}
\proof It follows from the proof of Lemma~\ref{s-lem3-1} that, for any $F\in\mathcal{F}\setminus\mathcal{F}_X$, we have $M=F\cup X,$ which implies that $F\in{M\choose k}.$  Let $\mathcal{A}^\prime=\{F\in {[n]\choose k}\mid X\subseteq F,\ |F\cap M|\geq t+1\}$ and $F^\prime$ be a fixed $k$-subset in $\mathcal{F}\setminus\mathcal{F}_X$. For any $F\in\mathcal{F}_X,$ since $|F\cap F^\prime|\geq t,$ $|F^\prime\cap X|\leq t-1$ and $M=F^\prime\cup X,$ we have $|F\cap M|\geq t+1,$ which implies that $\mathcal{F}_X\subseteq \mathcal{A}^\prime.$ Note that $\mathcal{A}^\prime\cup {M\choose k}$ is a $t$-intersecting family. By the maximality of $\mathcal{F},$ we have $\mathcal{F}=\mathcal{A}^\prime\cup {M\choose k}.$ $\qed$

By Remark~\ref{s-obs1}, if $\mathcal{F}$ is a maximal $t$-intersecting family satisfying the conditions in Lemma~\ref{s-prop3}, then $\mathcal{F}=\mathcal{H}_1(X,Y,M)$ for any $Y\in{M\choose k}$ with $X\subseteq Y$.
\begin{lemma}\label{s-prop3-1}
	Let $n,\ k,\ t,\ \mathcal{F}$ and $\mathcal{T}$ be as in Assumption~\ref{s-hyp1}. Suppose that $\tau_t(\mathcal{T})=t$, and $X$ is a $t$-subset of $[n]$ which is contained in each $T\in\mathcal{T}.$ Set $M=\cup_{T\in \mathcal{T}}T$ and $C=M\cup(\cup_{F\in\mathcal{F}\setminus \mathcal{F}_X}F)$.   Assume that $|M|=k$ and $|C|=c$. Then either $k+2\leq c\leq 2k-t$ or $c=n.$
	Moreover, the following hold.
	\begin{itemize}
		\item[{\rm(i)}] If $k+2\leq c\leq 2k-t,$ then $\mathcal{F}=\mathcal{H}_1(X,M,C).$
		\item[{\rm(ii)}]  If $c=n$, then $t\neq k-2$ and  $\mathcal{F}=\mathcal{H}_1(X,M,[n])$.
	\end{itemize}
\end{lemma}
\proof By the proof of Lemma~\ref{s-lem3-1}, for any $F\in\mathcal{F}\setminus\mathcal{F}_X,$ we have $|F\cap X|=t-1$ and $M\subseteq X\cup F$, which imply $|F\cap M|=k-1$ from $X\subseteq M.$ Choose $F_1\in \mathcal{F}\setminus\mathcal{F}_X.$ Then $|F_1\cup M|=k+1.$ If $c> k+1,$ then there exists $F_2\in \mathcal{F}\setminus\mathcal{F}_X$ such that $F_2\nsubseteq F_1\cup M,$ which implies that $F_2\cap (F_1\cup M)=F_2\cap M.$ Similarly, if $c>k+2,$ then there exists $F_3\in \mathcal{F}\setminus\mathcal{F}_X$ such that $F_3\nsubseteq F_1\cup F_2\cup M,$ which implies that $F_3\cap (F_1\cup F_2\cup M)=F_3\cap M.$ By mathematical induction, we can get $F_1,F_2,\ldots, F_{c-k}\in\mathcal{F}\setminus\mathcal{F}_X$ such that $F_i\cap(M\cup(\cup_{j=1}^{i-1}F_j))=F_i\cap M$ for any $i\in\{1,2,\ldots,c-k\}.$ If there exists $F^\prime\in\mathcal{F}$ such that $F^\prime\cap M=X$, then for any $i\in\{1,2,\ldots,c-k\}$, there exists $y_i\in F_i\setminus M$ such that $y_i\in F^\prime$ from $|F^\prime \cap F_i|\geq t$ and $|F^\prime \cap F_i\cap M|=t-1.$ Suppose $X=\{x_1,\ldots,x_t\}$. By the choice of $F_1,F_2,\ldots,F_{c-k},$ it is clear that $x_1,\ldots,x_t,y_1\ldots,y_{c-k}$ are  in $F^\prime$.

Suppose that $c\geq 2k-t+1.$ If there exists $F^\prime\in\mathcal{F}$ such that $F^\prime\cap M=X$, by above discussion, then $|F^\prime|>k$. That is impossible. Hence $|F^{\prime\prime}\cap M|\geq t+1$ for any $F^{\prime\prime}\in\mathcal{F}_X$. By the maximality of $\mathcal{F}$, it is easy to see that any $k$-subset $F^{\prime\prime\prime}$ of $[n]$ satisfying $|F^{\prime\prime\prime}\cap X|=t-1$ and $|F^{\prime\prime\prime}\cap M|=k-1$ is in $\mathcal{F}$. Then we have $C=[n]$ and $c=n.$ On the other hand, we have $c\geq k+2$, for otherwise we would have $c=k+1$, and $|T\cap F|\geq t$ for any $F\in\mathcal{F}$ and any $T\in{C\choose t+1}$ with $X\subseteq T$, which imply that $T\subseteq M$, a contradiction.

So far we have proved that either $k+2\leq c\leq 2k-t$ or $c=n.$ It remains to prove (i) and (ii).

(i)\quad Suppose that $k+2\leq c\leq 2k-t$.  Since  $|F\cap X|=t-1$ and $|F\cap M|=k-1$ for any $F\in\mathcal{F}\setminus\mathcal{F}_X,$ we have $\mathcal{F}\setminus\mathcal{F}_X\subseteq \mathcal{C}(X,M,C).$  For any $F^\prime\in\mathcal{F}_X,$ if $|F^\prime\cap M|\geq t+1$, then $F^\prime\in\mathcal{A}(X,M);$ if $F^\prime\cap M=X,$ then $|F^\prime\cap C|=c-k+t$ by above discussion, which implies that $F^\prime\in\mathcal{B}(X,M,C).$ Thus, $\mathcal{F}\subseteq \mathcal{H}_1(X,M,C).$ By the maximality of $\mathcal{F}$, we get $\mathcal{F}=\mathcal{H}_1(X,M,C)$.

(ii)\quad  Suppose that $c=n$.  Then $\mathcal{F}=\mathcal{A}(X,M)\cup\mathcal{C}(X,M,C)$ by the discussion in (i) and  maximality of $\mathcal{F}$. If $t=k-2$, then $\mathcal{F}=\mathcal{H}_1(X,M,[n])=\mathcal{H}_2(M)$, which implies that $\mathcal{T}={M\choose t+1}$ and $\tau_t(\mathcal{T})=t+1,$ a contradiction.  $\qed$
\begin{lemma}\label{s-prop4}
	Let $n,\ k,\ t,\ \mathcal{F}$ and $\mathcal{T}$ be as in Assumption~\ref{s-hyp1}. Suppose $\tau_t(\mathcal{T})=t+1$ and $\mathcal{T}={Z\choose t+1}$ for some $(t+2)$-subset $Z$ of $[n]$. Then $\mathcal{F}=\mathcal{H}_2(Z).$
\end{lemma}
\proof Since $\mathcal{T}={Z\choose t+1}$, we have $|F\cap Z|\geq t$ for any $F\in\mathcal{F}.$ If there exists $F^\prime\in\mathcal{F}$ such that $|F^\prime\cap Z|=t,$ then there exists a $T^\prime\in\mathcal{T}$ such that $|F^\prime\cap T^\prime|=t-1,$ a contradiction. Hence, $\mathcal{F}\subseteq \mathcal{H}_2(Z).$ Since $\mathcal{F}$ is maximal and $\mathcal{H}_2(Z)$ is $t$-intersecting, we have $\mathcal{F}=\mathcal{H}_2(Z)$. $\qed$

Now we give some equalities and inequalities for the sizes of Families I and II.
\begin{lemma}\label{s-lem1}
	 Suppose $c\in\{k+1,k+2,\ldots,2k-t,n\}.$ Then the following hold.
	\begin{align}
	h_1(t,k,c)=&{n-t\choose k-t}-{n-k\choose k-t}+{n-c\choose 2k-c-t}+t(c-k).\label{s-equ2}\\
	h_2(t+2)=&(t+2){n-t-2\choose k-t-1}+{n-t-2\choose k-t-2}.\label{s-equ4}
	\end{align}
\end{lemma}
\proof  Let $X,$ $M$, $C$, $\mathcal{A}(X,M)$, $\mathcal{B}(X,M,C)$ and $\mathcal{C}(X,M,C)$ be as in Family I. Then
\begin{align*}
|\mathcal{A}(X,M)|={n-t\choose k-t}-{n-k\choose k-t},\ |\mathcal{B}(X,M,C)|={n-c\choose 2k-c-t},\ |\mathcal{C}(X,M,C)|=t(c-k).
\end{align*}
Hence, (\ref{s-equ2}) holds.

Consider the family $\mathcal{H}_2(Z),$ where $Z$  is a $(t+2)$-subset of $[n]$. Observe that the number of $k$-subsets $F$ of $[n]$ satisfying $|F\cap Z|=t+1$ is $(t+2){n-t-2\choose k-t-1}$, and the number of $k$-subsets $F$ of $[n]$ satisfying $|F\cap Z|=t+2$ is ${n-t-2\choose k-t-2}$. Hence we have (\ref{s-equ4}) holds.    $\qed$

Let
\begin{align*}
f(n,k,t)=(k-t){n-t-1\choose k-t-1}-{k-t\choose 2}{n-t-2\choose k-t-2}.
\end{align*}
\begin{lemma}\label{s-lem3}
Let $1\leq t\leq k-2$ and $2k<n$. Then the following hold.
	\begin{itemize}
		\item[{\rm(i)}]  $h_1(t,k,k+1)>h_1(t,k,k+2)>\cdots>h_1(t,k,2k-t).$
		\item[{\rm(ii)}] $\min\{h_1(t,k,2k-t),\ h_1(t,k,n)\}\geq f(n,k,t).$
	\end{itemize}
\end{lemma}
\proof (i)\quad
For any $c\in\{k+1,k+2,\ldots,2k-t-1\}$, we have
$$
h_2(t,k,c)-h_2(t,k,c+1)={n-c-1\choose 2k-c-t}-t=\prod_{i=0}^{2k-c-t-1}\frac{n-c-1-i}{2k-c-t-i}-t.
$$
Observe that $\frac{n-c-1-i}{2k-c-t-i}>1$ for any $i\in\{0,1,\ldots,2k-c-t-1\},$ and $\frac{n-c-1-i}{2k-c-t-i}=n-2k+t>t$ when $i=2k-c-t-1$. Hence, we have $h_2(t,k,c)>h_2(t,k,c+1)$ for any $c\in\{k+1,k+2,\ldots,2k-t-1\}$, and (i) holds.
\
\medskip
\

(ii)\quad Let $X$ and $M$ be as in Family I. For any $i\in\{t,t+1,\ldots,k\}$, denote $
\mathcal{A}_i(X,M)=\left\{F\subseteq [n]\mid X\subseteq F,\ |F|=k,\ |F\cap M|=i\right\}$ and
$$
\mathcal{L}_i(X,M)=\left\{(I,F)\in{[n]\choose i}\times{[n]\choose k}\mid X\subseteq I\subseteq M,\ I\subseteq F \right\}.
$$
Double counting $|\mathcal{L}_i(X,M)|$,  we obtain
\begin{align*}
|\mathcal{L}_i(X,M)|=\sum_{j=i}^k\left|\mathcal{A}_j(X,M)\right|\cdot{j-t\choose i-t}={k-t\choose i-t}{n-i\choose k-i}.
\end{align*}
Since $\mathcal{A}(X,M)=\cup_{j=t+1}^k\mathcal{A}_j(X,M)$ and
$$
|\mathcal{L}_{t+1}(X,M)|=\sum_{j=t+1}^k\left|\mathcal{A}_j(X,M)\right|+\sum_{j=t+2}^k\left|\mathcal{A}_j(X,M)\right|\cdot\left(j-t-1\right),
$$
we obtain
\begin{align*}
|\mathcal{L}_{t+1}(X,M)|=&(k-t){n-t-1\choose k-t-1}\leq|\mathcal{A}|+\sum_{j=t+2}^k\left|\mathcal{A}_j(X,M)\right|\cdot {j-t\choose 2}\\
=&|\mathcal{A}|+|\mathcal{L}_{t+2}(X,M)|=|\mathcal{A}|+{k-t\choose 2}{n-t-2\choose k-t-2},
\end{align*}
which implies that $|\mathcal{A}(X,M)|\geq f(n,k,t)$. From the construction of $\mathcal{H}_1(X,M,C)$, we have (ii) holds.
$\qed$
\begin{lemma}\label{s-lem6}
	Let $1\leq t\leq k-2$ and ${t+2\choose 2}(k-t+1)^2+t\leq n$. Then the following hold.
	\begin{itemize}
		\item[{\rm(i)}] If $1\leq t\leq \frac{k}{2}-\frac{3}{2}$, then $h_2(t+2)<f(n,k,t)$.
		\item[{\rm(ii)}] If  $\frac{k}{2}-\frac{3}{2}< t\leq k-2$,  then $h_2(t+2)> f(n,k,t)$.
	\end{itemize}
\end{lemma}
\proof Note that
\begin{align}\label{s-equ3-4}
h_2(t+2)=(t+2){n-t-1\choose k-t-1}-(t+1){n-t-2\choose k-t-2}.
\end{align}
Let $f_2(n,k,t)=(f(n,k,t)-h_2(t+2))/{n-t-2\choose k-t-2}$. Then
$$
f_2(n,k,t)=\frac{(k-2t-2)(n-t-1)}{k-t-1}-{k-t\choose 2}+t+1
$$

\medskip

\noindent{\rm (i)}\quad From $1\leq t\leq \frac{k}{2}-\frac{3}{2}$ and ${t+2\choose 2}(k-t+1)^2+t\leq n$, we have
\begin{align*}
f_2(n,k,t)\geq \frac{(k-2t-2)(t+2)(t+1)(k-t+1)^2}{2(k-t-1)}-\frac{(k-t)(k-t-1)}{2}+t+\frac{t+1}{k-t-1}.
\end{align*}
Since $(k-2t-2)(t+2)(t+1)\geq (k-2t-2)+(t+2)=k-t$, we have $f_2(n,k,t)>0$ and (i) holds.

\medskip

\noindent{\rm (ii)}\quad If $t=\frac{k}{2}-1,$ then $f_2(n,k,t)=-t(t+1)/2<0.$ If $\frac{k}{2}-\frac{1}{2}\leq t\leq k-2,$ then $k-2t-2<0$ and
\begin{align*}
f_2(n,k,t)=&\frac{(k-2t-2)n+(t+1)^2}{k-t-1}-{k-t\choose 2}\\
\leq&\frac{(k-2t-2)(t+2)(t+1)(k-t+1)^2+2(t+1)^2}{2(k-t-1)}+\frac{(k-2t-2)t}{k-t-1}-{k-t\choose 2}<0.
\end{align*}
Therefore, we have (ii) holds.  $\qed$
\begin{lemma}\label{s-lem2}
Let $1\leq t\leq k-2$ and ${t+2\choose 2}(k-t+1)^2+t\leq n$. Then the following hold.
	\begin{itemize}
		\item[{\rm(i)}] Suppose that $1\leq t\leq k-3$. Then $h_1(t,k,2k-t-2)>h_1(t,k,n)\geq h_1(t,k,2k-t-1),$ and equality holds only if $t=1$.
		\item[{\rm(ii)}]  Suppose that $t=k-2$. Then $h_1(t,k,n)\geq h_1(t,k,k+1)$, and equality holds only if $t=1$.
	\end{itemize}
\end{lemma}
\proof  From Lemma~\ref{s-lem1}, we have
\begin{align*}
h_1(t,k,2k-t-2)-h_1(t,k,n)=&\frac{1}{2}(n-2k+t+2)(n-2k-t+1),\\
h_1(t,k,n)-h_1(t,k,2k-t-1)=&(t-1)(n-2k+t+1).
\end{align*}
Then (i) holds from
\begin{align*}
n-2k-t\geq{t+2\choose 2}(k-t+1)^2-2k\geq 2(t+2)(k-t+1)-2k=2(t+1)(k-t)+4>0.
\end{align*}
  If $t=k-2$,  then $h_1(t,k,n)-h_1(t,k,k+1)=(n-k-1)(t-1)$ and (ii) holds. $\qed$
\begin{lemma}\label{s-lem2-2}
Let $1\leq t\leq k-2$ and $2k<n.$ Then the following hold.
\begin{itemize}
	\item[{\rm(i)}] Suppose that $t=\frac{k}{2}-1$. Then $h_2(t+2)\geq h_1(t,k,k+2)$ and equality holds only if $t=1$; if $t=1$, or $t\geq 2$ and $n$ is sufficiently large, then $h_1(t,k,k+1)>h_2(t+2)$.
	\item[{\rm(ii)}] Suppose that $\frac{k}{2}-\frac{1}{2}\leq t\leq k-2$. Then $h_2(t+2)\geq h_1(t,k,k+1)$ and equality holds only if $(t,k)=(1,3)$.
\end{itemize}
\end{lemma}
\proof Since ${n-t\choose k-t}=\sum_{i=0}^{k-t-1}{n-k+i\choose k-t-1}+{n-k\choose k-t}$, we have
\begin{align}\label{s-equ3-3}
h_1(t,k,c)=\sum_{i=0}^{k-t-1}{n-k+i\choose k-t-1}+{n-c\choose 2k-c-t}+t(c-k)
\end{align}
for any $c\in\{k+1,k+2,\ldots,2k-t,n\}$. By (\ref{s-equ3-4}) and (\ref{s-equ3-3}), we have
\begin{align*}
&h_2(t+2)-h_1(t,k,k+1)\\=&(t+1)\left({n-t-1\choose k-t-1}-{n-t-2\choose k-t-2}\right)-\sum_{i=-1}^{k-t-2}{n-k+i\choose k-t-1}-t\\
=&(t+1){n-t-2\choose k-t-1}-\sum_{i=-1}^{k-t-2}{n-k+i\choose k-t-1}-t\\
=&(2t+1-k){n-t-2\choose k-t-1}+\sum_{i=-1}^{k-t-3}\left({n-t-2\choose k-t-1}-{n-k+i\choose k-t-1}\right)-t.
\end{align*}
If $\frac{k}{2}-\frac{1}{2}\leq t\leq k-3$, then $h_2(t+2)-h_1(t,k,k+1)\geq {n-t-2\choose k-t-1}-{n-k-1\choose k-t-1}-t>0.$ If $1=t=k-2,$ then $h_2(t+2)-h_1(t,k,k+1)=0$. If $2\leq t=k-2,$ then $h_2(t+2)-h_1(t,k,k+1)=(t-1)(n-t-2)+1-t>0$.

Suppose that $\frac{k}{2}-1=t$. If $t=1$, then $h_2(t+2)-h_1(t,k,k+1)=\frac{1}{2}(n-5)(-n+8)<0.$ When $t\geq2$, note that $h_2(t+2)-h_1(t,k,k+1)$ is a polynomial in $n$ with negative leading coefficient.  Then $h_2(t+2)-h_1(t,k,k+1)<0$ if $t\geq 2$ and $n$ is sufficiently large. By (\ref{s-equ3-4}) and (\ref{s-equ3-3}) again, we have
\begin{align*}
h_2(t+2)-h_1(t,k,k+2)=t{n-t-2\choose t+1}-\sum_{i=0}^{t-1}{n-2t-2+i\choose t+1}-{n-2t-4\choose t}-2t.
\end{align*}
If $t=1$, then $h_2(t+2)-h_1(t,k,k+2)=0$. If $t\geq 2,$ then
\begin{align*}
h_2(t+2)-h_1(t,k,k+2)>{n-t-2\choose t+1}-{n-2t-2\choose t+1}-{n-2t-4\choose t}-2t>0
\end{align*}
from ${n-t-3\choose t}>n-2t-2>2t$ and $
{n-t-2\choose t+1}={n-t-4\choose t+1}+{n-t-4\choose t}+{n-t-3\choose t}.
$ Hence, the desired result follows.   $\qed$

\section{Upper bounds for non-trivial $t$-intersecting families}
In \cite{Frankl-1978-1}, Frankl also proved that the size of any $t$-intersecting family $\mathcal{F}$ is no more than $c_k^\prime{n-\ell_1\choose k-\ell_1},$ where $\ell_1$ is the minimum size of a subset in the base of $\mathcal{F}$ and $c_k^\prime$ is a constant depending only on $k$. In this section, we give some specifical upper bounds on the sizes of the maximal non-trivial $t$-intersecting families. For any family $\mathcal{F}\subseteq{[n]\choose k}$ and any subset $S$ of $[n]$, define $
\mathcal{F}_S=\{F\in\mathcal{F}\mid S\subseteq F\}.$
\begin{lemma}\label{s-lem1-1}
	Let $\mathcal{F}\subseteq{[n]\choose k}$ be a $t$-intersecting family and $S$ an $s$-subset of $[n]$, where $t-1\leq s\leq k-1.$ If there exists $F^\prime\in \mathcal{F}$ such that $|S\cap F^\prime|=r\leq t-1,$ then for each $i\in\{1,2,\ldots,t-r\}$ there exists an $(s+i)$-subset $T_i$ with $S\subseteq T_i$ such that $|\mathcal{F}_S|\leq {k-r\choose i}|\mathcal{F}_{T_i}|$.
\end{lemma}
\proof For any $i\in\{1,2,\ldots,t-r\}$, let
\begin{align*}
\mathcal{H}_i=\{H\in{S\cup F^\prime \choose s+i}\mid S\subseteq H\}.
\end{align*}
 Observe that $|\mathcal{H}_i|={k-r\choose i}$. For any $F\in\mathcal{F}_S,$ since $\mathcal{F}$ is $t$-intersecting, we have $|F\cap F^\prime|\geq t$, implying that $|F\cap (S\cup F^\prime)|\geq s+t-r$ and there exists $H\in\mathcal{H}_i$ such that $H\subseteq F.$ Therefore $\mathcal{F}_S=\cup_{H\in\mathcal{H}_i}\mathcal{F}_H.$ Let $T_i$ be a subset in $\mathcal{H}_i$ such that $|\mathcal{F}_H|\leq|\mathcal{F}_{T_i}|$ for any $H\in\mathcal{H}_i.$ Thus $|\mathcal{F}_S|\leq {k-r\choose i}|\mathcal{F}_{T_i}|$ as desired.  $\qed$

Since $|\mathcal{F}_T|\leq {n-|T|\choose k-|T|}$ for any subset $T$ of $[n]$, we can obtain the following lemma.
\begin{lemma}\label{s-lem2-1}
	Let $\mathcal{F}\subseteq{[n]\choose k}$ be a $t$-intersecting family and $S$ an $s$-subset of $[n]$ with $t-1\leq s\leq k.$ If there exists $F^\prime\in \mathcal{F}$ such that $\dim(S\cap F^\prime)=r\leq t-1,$ then $|\mathcal{F}_S|\leq {k-r\choose t-r}{n-s-t+r\choose k-s-t+r}$.
\end{lemma}

The following lemma gives some upper bounds on the size of maximal non-trivial $t$-intersecting families $\mathcal{F}$ with $\tau_t(\mathcal{F})=t+1.$
\begin{lemma}\label{s-prop1}
	Let $n,\ k,\ t,\ \mathcal{F}$ and $\mathcal{T}$ be as in Assumption~\ref{s-hyp1}.  Then the following hold.
	\begin{itemize}
		\item[{\rm (i)}] If $|\mathcal{T}|=1$, then $|\mathcal{F}|\leq {n-t-1\choose k-t-1}+(t+1)(k-t)(k-t+1){n-t-2\choose k-t-2}$.
		\item[{\rm (ii)}] Suppose that $|\mathcal{T}|\geq 2$ and $	 \mathcal{T}=\left\{T\in{M\choose t+1}\mid X\subseteq T\right\}$  for some $t$-subset $X$ and $l$-subset $M$ of $[n]$ with $X\subseteq M.$ Then
		\begin{align}\label{s-equ9}
		|\mathcal{F}|\leq &(l-t){n-t-1\choose k-t-1}+(k-l+1)(k-t+1){n-t-2\choose k-t-2}+t{n-l\choose k-l+1}.
		\end{align}
		Moreover, if $l=t+2,$ then
		\begin{align}\label{s-equ10}
		|\mathcal{F}|\leq 2{n-t-1\choose k-t-1}+(k-1)(k-t+1){n-t-2\choose k-t-2}.
		\end{align}
		\item[{\rm (iii)}] If $|\mathcal{T}|\geq 2$ and  $\mathcal{T}={Z\choose t+1}$ for some  $(t+2)$-subset $Z$ of $[n]$, then $|\mathcal{F}|=h_2(t+2).$
	\end{itemize}
\end{lemma}
\proof (i) Let $T$ be the unique element in $\mathcal{T}$. Since $|T\cap F|\geq t$ for any $F\in\mathcal{F}$, we have
\begin{eqnarray}\label{s-equ7}
\mathcal{F}=\mathcal{F}_T\cup\left(\bigcup_{S\in{T\choose t}}(\mathcal{F}_S\setminus\mathcal{F}_T)\right).
\end{eqnarray}

We now give an upper bound on $|\mathcal{F}_S\setminus\mathcal{F}_T|$ for any fixed $S\in{T\choose t}.$ Since $\tau_t(\mathcal{F})=t+1,$ there exists an $F^\prime\in\mathcal{F}\setminus \mathcal{F}_S$ such that $|S\cap F^\prime|=t-1$ from $|F^\prime \cap T|\geq t.$ Then we have  $T=(F^\prime \cap T)\cup S$ and $T\subseteq F^\prime \cup S.$ For any $F\in \mathcal{F}_S\setminus\mathcal{F}_T$, notice that $(F\cap F^\prime)\cup S\subseteq F\cap (F^\prime\cup S).$ Since $|F\cap F^\prime|\geq t$ and $|F\cap F^\prime\cap S|\leq t-1$, we have $|F\cap(F^\prime\cup S)|\geq t+1.$ Hence there exists a $(t+1)$-subset $H$ such that $H\neq T$, $S\subseteq H\subseteq S\cup F^\prime$ and $H\subseteq F.$ Therefore, we have
\begin{eqnarray}\label{s-equ8}
\mathcal{F}_S\setminus\mathcal{F}_T=\bigcup_{S\subseteq H\subseteq S\cup F^\prime, \atop H\neq T, |H|=t+1} \mathcal{F}_H.
\end{eqnarray}

Consider any $(t+1)$-subset $H$ of $[n]$ satisfying $H\neq T$ and $S\subseteq H\subseteq S\cup F^\prime.$ Since $T$ is the unique $(t+1)$-subset of $[n]$ such that $|T\cap F|\geq t$ for any $F\in\mathcal{F}$, there exists $F^{\prime\prime}$ such that $|H\cap F^{\prime\prime}|<t,$ which implies that $|H\cap F^{\prime\prime}|=t-1$ from $|H\cap T|=|S|=t$ and $|T\cap F^{\prime\prime}|\geq t.$ From Lemma~\ref{s-lem2-1}, we have $|\mathcal{F}_{H}|\leq(k-t+1){n-t-2\choose k-t-2}.$ Observe $|\mathcal{F}_T|\leq{n-t-1\choose k-t-1}$ and
\begin{eqnarray*}
	\left|\left\{H\in{S\cup F^\prime\choose t+1}\mid S\subseteq H,\ H\neq T\right\}\right|=k-t.
\end{eqnarray*}
Therefore, from (\ref{s-equ7}) and (\ref{s-equ8}), we obtain
\begin{eqnarray*}
	|\mathcal{F}|\leq {n-t-1\choose k-t-1}+(t+1)(k-t)(k-t+1){n-t-2\choose k-t-2},
\end{eqnarray*}
as desired.
\
\medskip
\

(ii) We will obtain the upper bound of $|\mathcal{F}|$ by establishing upper bounds on $|\mathcal{F}_X|$ and $|\mathcal{F}\setminus\mathcal{F}_X|$. Since $\tau_t(\mathcal{F})=t+1,$ we have $|F\cap X|\geq t-1$ for any $F\in\mathcal{F}$, and there exists $F^\prime\in\mathcal{F}$ such that $|X\cap F^\prime|=t-1.$ From the proof of Lemma~\ref{s-lem3-1}, we have $X\subseteq M\subseteq X\cup F^\prime.$

For any $F\in\mathcal{F}_X,$ we have $|F\cap(X\cup F^\prime)|\geq t+1$ from $X\subseteq F$ and $|F\cap F^\prime|\geq t.$ So
\begin{eqnarray}\label{s-equ11}
\mathcal{F}_X=\left(\bigcup_{X\subseteq H_1,\ H_1\in{M\choose t+1}}\mathcal{F}_{H_1}\right)\cup\left(\bigcup_{X\subseteq H_2,\ H_2\in{X\cup F^\prime\choose t+1}\setminus{M\choose t+1}}\mathcal{F}_{H_2}\right).
\end{eqnarray}
Since $|\mathcal{F}_{H_1}|\leq {n-(t+1)\choose k-(t+1)}$ for any $H_1\in{M\choose t+1}$, we have $|\bigcup_{X\subseteq H_1,\ H_1\in{M\choose t+1}}\mathcal{F}_{H_1}|\leq(l-t){n-(t+1)\choose k-(t+1)}.$ For any $H_2\in{X\cup F^\prime\choose t+1}\setminus{M\choose t+1}$ with $X\subseteq H_2,$ since $H_2\notin \mathcal{T},$ there exists $F^{\prime\prime}\in\mathcal{F}$ such that $|H_2\cap F^{\prime\prime}|<t,$ which implies that $|H_2\cap F^{\prime\prime}|=t-1$ from $|F^{\prime\prime}\cap X|\geq t-1.$ It follows that $|\mathcal{F}_{H_2}|\leq (k-t+1){n-(t+1)-1\choose k-(t+1)-1}$ from Lemma~\ref{s-lem2-1}. Notice that
\begin{eqnarray*}
	\left|\left\{H_2\in{X\cup F^\prime\choose t+1}\setminus{M\choose t+1}\mid X\subseteq H_2\right\}\right|=k-l+1.
\end{eqnarray*}
Therefore, we have
\begin{eqnarray}\label{s-equ12}
|\mathcal{F}_X|\leq (l-t){n-t-1\choose k-t-1}+(k-l+1)(k-t+1){n-t-2\choose k-t-2}.
\end{eqnarray}

For any $F\in\mathcal{F}\setminus\mathcal{F}_X$ and any $T\in\mathcal{T}$, since $|F\cap X|=t-1$ and $X\nsubseteq F\cap T$, we have $T=(F\cap T)\cup X\subseteq F\cup X$. Then for any $F\in\mathcal{F}\setminus\mathcal{F}_X$ we have $M=\cup_{T\in\mathcal{T}}T\subseteq F\cup X,$ which implies that $|M\cap F|=l-1$. Hence, $\mathcal{F}\setminus\mathcal{F}_X\subseteq\{F\in{[n]\choose k}\mid |F\cap M|=l-1,\ X\nsubseteq F\},$ and
\begin{align}
|\mathcal{F}\setminus\mathcal{F}_X|\leq t{n-l\choose k-l+1}. \label{s-equ13}
\end{align}

Combining  (\ref{s-equ12}) and (\ref{s-equ13}), we obtain (\ref{s-equ9}).

Now let us consider the case when $l=t+2$. From the discussion above, we have $|M\cap F|=l-1=t+1$ for any $F\in\mathcal{F}\setminus\mathcal{F}_X,$ which implies that $$\mathcal{F}\setminus\mathcal{F}_X\subseteq \bigcup_{X\nsubseteq L,\ L\in{M\choose t+1}} \mathcal{F}_L.$$
For any $L\in{M\choose t+1}$ with $X\nsubseteq L$, since $L\notin\mathcal{T}$ and $|F\cap M|\geq t$ for any $F\in\mathcal{F}$, there exists $F^\prime\in\mathcal{F}$ such that $|F^\prime\cap L|=t-1.$ Then $|\mathcal{F}_L|\leq (k-t+1){n-t-2\choose k-t-2}$ from Lemma~\ref{s-lem2-1}. Since the number of $(t+1)$-subsets $L$ of $M$ with $X\nsubseteq L$ is equal to $t$, we have
\begin{eqnarray}\label{s-equ14}
|\mathcal{F}\setminus\mathcal{F}_X|\leq t(k-t+1){n-t-2\choose k-t-2}.
\end{eqnarray}
Combining (\ref{s-equ12}) and (\ref{s-equ14}), we obtain (\ref{s-equ10}).
\
\medskip
\

(iii) By Lemma~\ref{s-prop4}, the desired result follows.  $\qed$
\begin{lemma}\label{s-prop6}
Let $n$, $k$ and $t$ be integers with $1\leq t\leq k-2$ and ${t+2\choose 2}(k-t+1)^2+t\leq n$, and let $\mathcal{F}\subseteq {[n]\choose k}$ be a maximal $t$-intersecting family with $t+2\leq \tau_t(\mathcal{F})=m\leq k$. Then
\begin{align}\label{s-equ15-1}
|\mathcal{F}|\leq k^{m-t-2}(k-t+1)^2{m\choose t}{n-m\choose k-m}.
\end{align}	
Moreover, we have
\begin{align}\label{s-equ16-1}
|\mathcal{F}|\leq (k-t+1)^2{t+2\choose 2}{n-t-2\choose k-t-2}.
\end{align}
\end{lemma}
\proof Let $T$ be an $m$-subset of $[n]$ which satisfies $|T\cap F|\geq t$ for any $F\in\mathcal{F}$. Then $\mathcal{F}=\cup_{H\in{T\choose t}}\mathcal{F}_H$ and there exists $H_1\in{T\choose t}$ such that $|\mathcal{F}|\leq {m\choose t}|\mathcal{F}_{H_1}|$. If $m\geq t+3,$ using Lemma~\ref{s-lem1-1} repeatedly, then there exist $H_2\in{[n]\choose t+1}$, $H_3\in{[n]\choose t+2}$,\ldots,$H_{m-t-1}\in{[n]\choose m-2}$ such that $H_i\subseteq H_{i+1}$ and $|\mathcal{F}_{H_i}|\leq k|\mathcal{F}_{H_{i+1}}|$ for each $i\in\{1,2,\ldots,m-t-2\}$. Thus there exists $H^\prime\in{[n]\choose m-2}$ such that
\begin{align*}
|\mathcal{F}|\leq{m\choose t}k^{m-t-2}|\mathcal{F}_{H^\prime}|.
\end{align*}
Since $\tau_t(\mathcal{F})>m-2$, we have $\mathcal{F}\setminus\mathcal{F}_{H^\prime}\neq\emptyset$ and $|F\cap H^\prime|\leq t-1$ for any $F\in\mathcal{F}\setminus\mathcal{F}_{H^\prime}.$

\medskip

\noindent\textbf{Case 1.} $|F\cap H^\prime|\leq t-2$ for all $F\in\mathcal{F}\setminus\mathcal{F}_{H^\prime}$.

In this case, we have $t\geq 2$. For $s\in\{0,1,\ldots,t-2\}$, let
\begin{align*}
g(s)={k-s\choose t-s}{n-m+2-t+s\choose k-m+2-t+s}.
\end{align*}
Since $t+2\leq m$, $s\leq t-2$ and ${t+2\choose 2}(k-t+1)^2+t\leq n$, we have
\begin{align*}
&(t-s)(n-m+3-t+s)-(k-s)(k-m+3-t+s)\\
=&(k-t)m+(t-s)(n-k)-(k-t)(k+3-t+s)\\
>&(k-t)(t+2)+n-k-(k-t)(k+1)\\
=&n-t-(k-t)^2>0,
\end{align*}
which implies that
\begin{align*}
\frac{g(s+1)}{g(s)}=\frac{(t-s)(n-m+3-t+s)}{(k-s)(k-m+3-t+s)}>1
\end{align*}
for $s\in\{0,1,\ldots,t-3\}$. That is the function $g(s)$ is increasing as $s\in\{0,1,\ldots,t-2\}$ increases.

Let $F_1$ be a fixed $k$-subset in $\mathcal{F}\setminus\mathcal{F}_{H^\prime}$. Assume that $|F_1\cap H^\prime|=s_1$. Observe that $0\leq s_1\leq t-2.$ By Lemma~\ref{s-lem2-1}, we have
$|\mathcal{F}_{H^\prime}|\leq g(s_1)\leq g(t-2)$, which implies that
\begin{align}\label{s-equ15}
|\mathcal{F}|\leq {m\choose t}k^{m-t-2}g(t-2)=k^{m-t-2}{m\choose t}{k-t+2\choose 2}{n-m\choose k-m}.
\end{align}

\medskip

\noindent\textbf{Case 2.} There exists $F_2\in\mathcal{F}\setminus\mathcal{F}_{H^\prime}$ such that $|F_2\cap H^\prime|=t-1$.

By Lemma~\ref{s-lem1-1}, there exists an $(m-1)$-subset $H^{\prime\prime}$ such that $|\mathcal{F}_{H^\prime}|\leq (k-t+1)|\mathcal{F}_{H^{\prime\prime}}|$. Hence, we have $|\mathcal{F}|\leq{m\choose t}k^{m-t-2}(k-t+1)|\mathcal{F}_{H^{\prime\prime}}|$. Since $\tau_t(\mathcal{F})>m-1$, there exists $F_3\in\mathcal{F}$ such that $|F_3\cap H^{\prime\prime}|\leq t-1$.

If $|F_3\cap H^{\prime\prime}|=t-1,$ then there exists an $m$-subset $H^{\prime\prime\prime}$ of $[n]$ with $H^{\prime\prime}\subseteq H^{\prime\prime\prime}$ such that $|\mathcal{F}_{H^{\prime\prime}}|\leq (k-t+1)|\mathcal{F}_{H^{\prime\prime\prime}}|$. Since $|\mathcal{F}_{H^{\prime\prime\prime}}|\leq{n-m\choose k-m}$, we have
\begin{align}\label{s-equ16}
|\mathcal{F}|\leq k^{m-t-2}(k-t+1)^2{m\choose t}{n-m\choose k-m}.
\end{align}

Suppose that $t\geq 2$ and $|F_3\cap H^{\prime\prime}|=s_2\leq t-2$. By Lemma~\ref{s-lem2-1}, we have
\begin{align*}
|\mathcal{F}_{H^{\prime\prime}}|\leq {k-s_2\choose t-s_2}{n-m+1-t+s_2\choose k-m+1-t+s_2}.
\end{align*}
Similar to Case 1, it is straightforward to verify that the function ${k-s\choose t-s}{n-m+1-t+s\choose k-m+1-t+s}$ is increasing as $s\in\{0,1,\ldots,t-2\}$ increases. Hence
\begin{align}\label{s-equ17}
|\mathcal{F}|\leq k^{m-t-2}(k-t+1){m\choose t}{k-t+2\choose 2}{n-m-1\choose k-m-1}.
\end{align}

If $t=1$, then (\ref{s-equ15-1}) holds from (\ref{s-equ16}). If $t\geq 2$, from ${t+2\choose 2}(k-t+1)^2+t\leq n$, it is straightforward to verify that
\begin{align*}
(k-t+1)^2{n-m\choose k-m}\geq\max\left\{{k-t+2\choose 2}{n-m\choose k-m},\ (k-t+1){k-t+2\choose 2}{n-m-1\choose k-m-1}\right\},
\end{align*}
which together with (\ref{s-equ15}), (\ref{s-equ16}) and (\ref{s-equ17}) yields that (\ref{s-equ15-1}) holds.

Let $p(x)=(x-t+1)(n-x)-k(x+1)(k-x)$ for each $x\in\{t+1,t+2,\ldots,k\}$. Observe that
\begin{align*}
p(t+1)=&2(n-t-1)-k(t+2)(k-t-1)\\
\geq& (t+2)(t+1)(k-t+1)^2-k(t+2)(k-t+1)+2k(t+2)-2>0,
\end{align*}
and
\begin{align*}
p(x+1)-p(x)=&n-k^2+2k+t-2+(2k-2)x\\
\geq&(k-t+1)^2+t-k^2+2k+t-2+(2k-2)(t+1)>0
\end{align*}
for each $x\in\{t+1,t+2,\ldots,k-1\}.$ Hence $p(x)>0$ for any $x\in\{t+1,t+2,\ldots,k\}$. Let
\begin{align*}
q(y)=k^{y-t-2}{y\choose t}{n-y\choose k-y}
\end{align*}
for each $y\in\{t+2,t+3,\ldots,k\}$. From ${t+2\choose 2}(k-t+1)^2+t\leq n$ and $p(x)>0$ for any $x\in\{t+1,t+2,\ldots,k\}$, we have
\begin{align*}
\frac{q(y)}{q(y+1)}=\frac{(y-t+1)(n-y)}{k(y+1)(k-y)}>1
\end{align*}
for any $y\in\{t+2,t+2,\ldots,k-1\}$. That is, the function $q(y)$ is decreasing as $y\in\{t+2,t+3,\ldots,k\}$ increases. This together with (\ref{s-equ15-1}) yields (\ref{s-equ16-1}) holds.     $\qed$

\section{The proof of Theorem~\ref{s-main-1}}
Let $\mathcal{F}$ be any maximal non-trivial $t$-intersecting family which is not given in  Theorem~\ref{s-main-1}. Suppose $f_3(n,k,t)=\left(f(n,k,t)-|\mathcal{F}|\right)/{n-t-2\choose k-t-2}.$ It suffices to prove that $f(n,k,t)>|\mathcal{F}|$ or $f_3(n,k,t)>0$.

\noindent\textbf{Case 1.}\quad $\tau_t(\mathcal{F})=t+1.$
\
\medskip
\

Let $\mathcal{T}$ be the set of all $(t+1)$-subsets $T$ of $[n]$ which satisfies $|T\cap F|\geq t$ for any $F\in\mathcal{F}$.

Suppose that $|\mathcal{T}|=1$. Since $n\geq{t+2\choose 2}(k-t+1)^2+t$, by Lemma~\ref{s-prop1} (i), we have
$$
f_3(n,k,t)\geq n-t-1-{k-t\choose 2}-(t+1)(k-t)(k-t+1)>0.
$$

Suppose that $|\mathcal{T}|\geq 2$ and $\tau_t(\mathcal{T})=t$. Assume that $l=t+2$. By Lemma~\ref{s-prop1} (ii), we have
\begin{align*}
f_3(n,k,t)\geq \frac{(k-t-2)(n-t-1)}{k-t-1}-{k-t\choose 2}-(k-1)(k-t+1).
\end{align*}
 Observe that $n\geq{t+2\choose 2}(k-t+1)^2+t$. If $k=t+3,$ then $f_3(n,k,t)\geq 4(k^2-4k+2)+0.5>0.$ If $k=t+4,$ then $f_3(n,k,t)\geq \frac{5}{3}(5k^2-28k+29)>0.$ If $k>t+4,$ then $n-t-1>(t+2)(k-t+1)(k-t-1)$ and
\begin{align*}
f_3(n,k,t)>&(k-t-2)(t+2)(k-t+1)-(k-1)(k-t+1)-(k-t)(k-t-1)/2\\
=&(k-t-4)(k-t+1)t+(k-3)(k-t+1)-(k-t)(k-t-1)/2\\
=&(k-t-4)(k-t+1)t+(k-t)(k+t-5)/2+(k-3)>0.
\end{align*}
Assume that $t+2<l<k$. Since $\max\left\{{t+2\choose 2},\frac{k-t+2}{2}\right\}\geq \frac{1}{k-t}{t+2\choose 2}+\frac{(k-t-1)(k-t+2)}{2(k-t)}$ and  ${n-l\choose k-l+1}<{n-t-2\choose k-t-2}$, By Lemma~\ref{s-prop1} (ii), then
\begin{align*}
f_3(n,k,t)>& \frac{(k-l)(n-t-1)}{k-t-1}-{k-t\choose 2}-(k-l+1)(k-t+1)-t\\
\geq& \frac{n-t-1}{k-t-1}-{k-t\choose 2}-2(k-t+1)-t\\
\geq & \frac{(k-t+1)^2}{(k-t)(k-t-1)}{t+2\choose 2}-t-1+\frac{(k-t+2)(k-t+1)^2}{2(k-t)}-{k-t\choose 2}-2(k-t+1)\\
>&1+\frac{(k-t+2)(k-t+1)}{2}-{k-t\choose 2}-2(k-t+1)\geq0.
\end{align*}

Suppose $|\mathcal{T}|\geq 2$ and $\tau_t(\mathcal{T})=t+1$. By Lemmas~\ref{s-prop4} and \ref{s-lem6}, we have $f(n,k,t)>|\mathcal{F}|$ if $1\leq t\leq \frac{k}{2}-\frac{3}{2}$.
\
\medskip
\

\noindent\textbf{Case 2.}\quad $t+2\leq\tau_t(\mathcal{F})\leq k.$
\
\medskip
\

Since
\begin{align*}
\max\left\{{t+2\choose 2},\ \frac{k-t+2}{2}\right\}\geq \frac{k-t-1}{k-t}{t+2\choose 2}+\frac{1}{k-t}\cdot\frac{k-t+2}{2},
\end{align*}
by Lemma~\ref{s-prop6}, we have
\begin{align*}
f_3(n,k,t)=&\frac{(k-t)(n-t-1)}{k-t-1}-{k-t\choose 2}-(k-t+1)^2{t+2\choose 2}\\
\geq&\frac{(k-t+1)^2(k-t+2)-2(k-t)}{2(k-t-1)}-{k-t\choose 2}>0.
\end{align*}
Hence the desired result follows.  $\qed$

\section*{Acknowledgement}

This research is supported by NSFC (11671043) and NSF of Hebei Province(A2019205092).

%\end{CJK*}

\end{document}